\documentclass{article}

\usepackage{amssymb}

\newcommand{\alsxtilde}{\tilde{\alpha} _{S} (x)}

\newcommand{\altilde}{\tilde{\alpha}(x)}

\newcommand{\alx}{\alpha (x)}

\newcommand{\demo}{{\bf Proof}}

\newcommand{\dphi}{d \varphi}

\newcommand{\entero}{\mathbb{Z}}

\newcommand{\ex}{ex}

\newcommand{\fibra}{S(x) = \txSigma / \xhx}

\newcommand{\fibvec}{\pi : S \to \Sigma}

\newcommand{\fibvecsx}{\pi : S \to X}

\newcommand{\fs}{\cal F _{\cal S}}

\newcommand{\fstilde}{\tilde{\fs}}

\newcommand{\ham}{H:M \to \numreal}

\newcommand{\htop}{h_{top}}

\newcommand{\LsSigma}{\Lambda (S)}

\newcommand{\neps}{N_{\epsilon}}

\newcommand{\nepsprima}{\neps '}

\newcommand{\numreal}{\mathbb{R}}

\newcommand{\oms}{\omega_{S}}

\newcommand{\omSigma}{\omega _{\Sigma}}

\newcommand{\phit}{\varphi _{t}}

\newcommand{\phitast}{\phi _{t} ^{\ast}}

\newcommand{\px}{p_{x} : \txSigma \to \sx}

\newcommand{\sdos}{S^{2}}

\newcommand{\sx}{S(x)}

\newcommand{\tuno}{\tau _{1}}

\newcommand{\tdos}{\tau _{2}}

\newcommand{\txSigma}{T_{x} \Sigma}

\newcommand{\vuno}{v_{1}}

\newcommand{\vdos}{v_{2}}

\newcommand{\xhx}{X_{H} (x)}

\newcommand{\zetaeta}{\zeta, \eta}

\newtheorem{prop}{\bf{Proposition}}

\newtheorem{lema}{\bf{Lemma}}

\newtheorem{teo}{\bf{Theorem}}

\title{On the topological entropy of an optical Hamiltonian flow}

\author{C\'esar J. Niche \\ Mathematics Department \\ University Of California at Santa Cruz \\ Santa Cruz CA 95064 \\ USA \\ E-mail: cniche@math.ucsc.edu }

\begin{document}

\noindent

\maketitle

\begin{abstract}

In this article, we prove two formulas for the topological entropy of an ${\cal F}$-optical Hamiltonian flow induced by $H \in C^{\infty} (M, \numreal)$, where ${\cal F}$ is a Lagrangian distribution. In these formulas, we calculate the topological entropy as the exponential growth rate of the average of the determinant of the differential of the flow, restricted to the Lagrangian distribution or to a proper modification of it.

\end{abstract}

\section{Introduction}
\label{intro}

Bounds or formulas for calculating the topological entropy of a dynamical system 
are very important tools when studying diffeomorphisms or flows on compact 
manifolds. Many results exist relating this conjugacy invariant to the growth 
rate of volumes of submanifolds or to the growth rate of the expansion on the 
tangent bundle, for instance, Cogswell \cite{cogswell}, Kozlovski \cite{k}, Newhouse 
\cite{newhouse2}, Przytycki \cite{prz} and Yomdin \cite{yomdin}. \\

For a closed Riemannian manifold $(M, g)$, where the metric $g$ is $C^{\infty}$, 
Ma\~n\'e \cite{manhe} proved the following results for the topological entropy of 
the geodesic flow $\phit:SM\to SM$, in terms of $n_{T}(x,y)$, the number of geodesics of length 
less or equal than $T$ between two points  $x$ and $y$,  of the 
expansion

\begin{displaymath}
ex \, (\dphi_{T})_{\theta} = \max _{S \subset T_{\theta}SM} \vert det \, 
(\dphi_{T})_{\theta} \vert _{S} \vert
\end{displaymath}

and of $V(\theta) = Ker (d\pi)_{\theta}$, where $\pi: SM \to M$ is the usual 
projection.

\begin{teo} [Theorem 1.1, Ma\~n\'e \cite{manhe}]

\label{teomanhe1}
\begin{displaymath}
\htop = \lim _{T \to \infty} \frac{1}{T} \log \int _{M \times M} n_{T} (x,y) \, 
dx \, dy = \lim _{T \to \infty} \frac{1}{T} \log \int _{SM} \ex 
(\dphi_{T})_{\theta} \, d\theta.
\end{displaymath}

\end{teo}

\begin{teo} [Teorema 1.4, Ma\~n\'e \cite{manhe}]  

\label{teomanhe2}
\begin{displaymath}
\htop = \lim _{T \to \infty} \frac{1}{T} \log \int _{SM} \vert det \, (\dphi 
_{T})_{\theta} \vert _{V(\theta)} \vert \, d\theta = \lim _{T \to \infty} 
\frac{1}{T} \log \int _{M} vol \, \varphi_{T} (S_{x}M) \, dx.
\end{displaymath}

\end{teo}

In spite of the fact that these results are Riemannian in nature, the second 
equality in Theorem \ref{teomanhe1} and the first in Theorem \ref{teomanhe2} are 
proved by symplectic techniques, relying on a canonical symplectic structure on 
$TM$, which allows us to treat the geodesic flow as a Hamiltonian one. For a 
clear exposition of these and related results, see Paternain \cite{pat}.

There are three key facts on the proof of Ma\~n\'e 's results. They are: the 
twist property of the Lagrangian vertical subspace $V(\theta)$, which says that 
if the intersection of $V (\phit (\theta))$ and $(d \phit) _{\theta} (V(\theta))$ 
is nontrivial, a slight perturbation in $t$ makes it trivial; Przytycki 's 
inequality \cite{prz}, which gives a bound for the topological entropy of a 
$C^{2}$ flow in terms of the exponential growth rate of the average of the 
expansion on $T(SM)$ and a very nice change of variables and some auxiliary 
lemmas, which lead to an inequality necessary to bound from above the average of 
the expansion. \\

It is natural then, to try to prove similar results for certain Hamiltonian 
flows, taking into account the importance of Lagrangian distributions in 
symplectic geometry. To proceed with such a generalization, the vertical 
distribution $V(\theta)$, Przytycki's inequality and the technical lemmas, should be properly replaced within the Hamiltonian context. \\   

In order to recover a twist property, we need to introduce the concept of optical 
Hamiltonian. We closely follow Bialy and Polterovich \cite{bp}, \cite{bp2}.

Given a symplectic vector space $(E^{2n}, \omega)$, let $\Lambda (E)$ be the 
family of its Lagrangian subspaces, which is a manifold diffeomorphic to  
$U(n)/O(n)$. For $\lambda \in \Lambda$, the tangent space $T_{\lambda}\Lambda$ can 
be canonically identified with $\sdos(\lambda)$, the vector space of bilinear  
symmetric forms on $\lambda$, in the following way: given a curve $\lambda(t) 
\subset \Lambda$ with $\lambda(0) = \lambda$, it is described by 

\begin{equation}
\lambda(t) = S(t) \lambda, \qquad S(0) \, = \, Id, \, S(t) \in Sp(E, \omega)
\end{equation}

where $Sp(E, \omega)$ is the group of symplectic maps in $E$. Thus, to the vector $\dot{\lambda} (0) \in T_{\lambda}\Lambda$ we can associate 
the symmetric form

\begin{equation}
\label{eq:df-forma}
(\zetaeta) \mapsto \omega(\zeta, \dot{S}(0) \eta), \qquad \zeta, \eta \in 
\lambda. \\
\end{equation}

For a fixed Lagrangian space $\alpha \in \Lambda$, we note as 

\begin{equation}
\label{eq:estrato}
\Lambda _{\alpha} = \bigcup _{1 \leq k \leq n} \Lambda _{\alpha} ^{(k)}
\end{equation}

the stratified manifold where  $\Lambda _{\alpha} ^{(k)}$ is the family of 
Lagrangian subspaces whose intersection with  $\alpha$ is $k$-dimensional. So, if 
$\lambda \in \Lambda _{\alpha} ^{(k)}$, we can identify  $T_{\lambda}\Lambda / 
T_{\lambda}\Lambda _{\alpha} ^{(k)}$ with $S^{2} (\alpha \cap \lambda)$. Then, for 
$\lambda \in \Lambda _{\alpha} ^{(k)}$, we say that a vector in  
$T_{\lambda}\Lambda / T_{\lambda}\Lambda _{\alpha} ^{(k)}$ is $\alpha$-positive 
if the bilineal form (\ref{eq:df-forma}) associated to it is positive definite.   
A vector in $T_{\lambda}\Lambda$ is $\alpha$-positive if its image under the 
projection to $T_{\lambda}\Lambda / T_{\lambda}\Lambda _{\alpha} ^{(k)}$ is 
$\alpha$-positive. \\

Generally, to a symplectic fibre bundle $\pi : E \to M$, we can associate 
its Lagrange-Grassman fibre bundle $\Theta: {\cal A} \to M$, where each fibre is 
${\cal A} _{x} = \Lambda (T_{x}M)$. Given a Lagrangian distribution $\cal F$, that is, a section ${\cal F}: 
M \to {\cal A}$, we see that the subbundle $\cal A _{\cal F}$ of Lagrangian 
subspaces not tranversal to $\cal F$, has fibres diffeomorphic to $\Lambda 
_{\alpha}$ in (\ref{eq:estrato}). Then, a vector in $T _{\lambda} \cal A$ is 
$\cal F$-positive for $\lambda \in \cal A _{\cal F}$, if its image under the projection $pr: T_{\lambda}{\cal A} \to T_{\lambda}\Lambda / T_{\lambda}\Lambda _{\alpha} ^{(k)}$ is $\cal F$-positive. A 
differentiable curve in $\cal A$, is $\cal F$-positive, if its tangent vectors 
are $\cal F$-positive where they intersect $\cal A _{\cal F}$ .

The following proposition gives us the proper notion of twist property. 

\begin{prop} [Prop. 1.6, Bialy y Polterovich \cite{bp}]
\label{propbp}
Every ${\cal F}$-positive tangent vector to $\cal A$, is transversal to the 
stratified manifold $\cal A_{\cal F}$.
\end{prop}

We now define optical Hamiltonians. We consider the curve of 
Lagrangian subspaces $(\phit)_{*} (\alpha (\varphi _{-t} (x))$ on $\pi: {\cal A} 
\to M$, where the base $M$ is a symplectic manifold endowed with a Lagrangian 
distribution $\cal F$, and $\phit$ is the flow induced by $\ham$. Then $H$ 
is $\cal F$-optical if for every $x \in M, t \in \numreal$, 

\begin{displaymath}
\frac{d}{dt} \, (\phit)_{*} \, (\alpha (\varphi _{-t} (x))) \, \vert _{t = 0}
\end{displaymath}

is an $\alx$-positive vector. A Hamiltonian flow  is optical, if it is generated 
by an optical Hamiltonian. \\

We now state the two main results of this article. Let $(M, \omega)$ be a 
symplectic manifold, endowed with a Lagrangian distribution $\cal F$ of 
Lagrangian subspaces $\alx \subset T_{x}M$. Let $H \in C^{\infty}(M, \numreal)$
be an $\cal F$-optical Hamiltonian with induced flow $\phit$. Let $\Sigma = H 
^{-1} (e)$ be a compact energy level, for a regular value $e$.

\begin{teo}
\label{teoa}
For $\neps = H^{-1} (e - \epsilon, e + \epsilon)$, then

\begin{displaymath}
\htop (\phit) = \lim _{\epsilon \to 0} \lim _{t \to \infty} 
\frac{1}{t} \log \int _{\neps} \vert det \, (\dphi _{t})_{x} \vert _{\alx} \vert 
\,dx.
\end{displaymath}

\end{teo}

\begin{teo}
\label{teob} 
Let us assume that there exists a continuous invariant distribution of hyperplanes 
$T(x)$ transversal to the the Hamiltonian vector field on $\Sigma$. Then 

\begin{displaymath}
\htop (\phit) = \lim _{t \to \infty} \frac{1}{t} \log \int 
_{\Sigma} \vert det \, (\dphi _{t})_{x} \vert _{\altilde} \vert \,dx
\end{displaymath}

where $\altilde = (\alx \cap \txSigma) + \langle \xhx \rangle$.
\end{teo}

The proof of these results is based upon Proposition \ref{prop-clave-prev}, which 
is a slightly more general version of the bound for the average of the expansion  already mentioned. We state this result now. Given a 
symplectic linear cocycle, that is, a $C^{\infty}$ flow on a compact manifold $X$ 
and a family of differentiable maps $\phitast (x): \sx \to S(\phi _{t} (x))$, such 
that $\phi _{t} \circ \pi = \pi \circ \phi _{t}^{*}$ and that for all $x \in X, t 
\in \numreal$, $\phitast (x)$ is a linear symplectic isomorphism, then we say the 
cocycle is $\cal F$-optical with respect to a Lagrangian distribution $\cal 
F$, if it is in the sense of the previous paragraphs.

\begin{prop}
\label{prop-clave-prev}
If $\phitast$ is $\cal F$-optical, then there exists a constant $C > 0$ such that 
for all $t \in \numreal$

\begin{displaymath}
\int_{X} \vert det \, \phitast(x) \vert_{\alx} \vert \, dx \geq C \, \int _{X}  
ex \, \phitast(x) \, dx. \\
\end{displaymath}
\end{prop}

These Theorems are interesting given the generality of the opticity condition. 
For example, on $(T^{\ast}M, \omega _{0})$ where $\omega _{0}$ is the canonical symplectic form, a Hamiltonian $H \in C^{\infty} 
(T^{\ast}M, \numreal)$, is $\cal F$-optical with respect to the distribution ${\cal 
F} = \{ dq = 0 \}$ of spaces tangent to the fibres iif  $H$ is convex 
on each fibre, that is,  iif $H _{pp} > 0$.  As a result of this, the usual 
Hamiltonians of classical mechanics, are  optical for this distribution. 

Another important case where these results can be applied, is that of Anosov magnetic 
flows induced by twisted symplectic structures. See for instance, Burns and 
Paternain \cite{burnspat}. \\

This article is organized as follows. In  \S \ref{previos} we prove Proposition 
\ref{prop-clave-prev} and then in \S \ref{pruebateoa} and \S \ref{pruebateob}, we prove Theorems \ref{teoa} and  \ref{teob}. \\

{\bf Acknowledgements} We deeply thank Gabriel Paternain for his advice. The author acknowledges partial financial support 
from PEDECIBA - Area Matem\'a\-tica, Uruguay and from Instituto de Matem\'atica y Estad\'{\i}stica Rafael Laguardia (IMERL), Facultad de Ingenier\'{\i}a, Universidad de la Rep\'ublica, Montevideo, Uruguay. An important part of this article was 
written during a stay at Centro de Investigaci\'on en Matem\'atica (CIMAT), 
Guanajuato, M\'exico. This stay was supported by Comisi\'on Sectorial de 
Investigaci\'on Cient\'{\i}fica (CSIC) from Universidad de la Rep\'ublica, 
Montevideo, Uruguay and CIMAT. To Renato Iturriaga (CIMAT) and to these 
institutions, we are deeply grateful.

\section{Proof of Proposition \ref{prop-clave-prev}}

\label{previos}

Before proving Proposition \ref{prop-clave-prev}, we recall some definitions and results needed. \\

In order to measure angles and volumes on a symplectic manifold $(M, \omega)$, 
we introduce a compatible Riemannian metric. Given a symplectic vector bundle 
$\pi: E \to M$, an almost complex structure is $J: E \to E$, such that $J^{2} = 
- Id _{E}$. It is said compatible with the symplectic form $\omega$, if on each fiber 
$E_{q}$

\begin{eqnarray}
\omega (v,w) = \omega (Jv, Jw), \quad \forall v, w \in E_{q} \nonumber \\
\omega(v, Jv) > 0, \quad \forall v \in E_{q}.
\end{eqnarray}

Thus, we can define a Riemannian metric on $E$ by $g_{J} (v,w) = \omega (v, 
Jw)$. For a proof of the existence and properties of such an almost complex structure, see Mc Duff y Salamon 
\cite{mcd-s}. \\

Given a linear map $L: E \to F$, where $E$ y $F$ are finite dimensional  
Hilbert spaces, we define the determinant of $L$ as follows. For an 
orthonormal base of $E$,  $\{ \vuno, \vdos, \ldots, v_{n} \}$, we consider the 
matrix $a_{ij} = \langle L(v_{i}), L(v_{j}) \rangle$. Then, we define 

\begin{displaymath}
\vert det \, L \vert = \sqrt{det \, A}.
\end{displaymath}

We also define the expansion of $L$ as 

\begin{displaymath}
ex \, L = \max _{S \subset E} \vert det L \vert _{S} \vert
\end{displaymath}

where $S$ is a subspace in $E$.

Finally, for $E_{1}, E_{2}$ subspaces of $E$, with $dim \, E_{i} = \frac{1}{2} \, 
dim \, E$, we define the angle between them as 

\begin{displaymath}
ang (E_{1}, E_{2}) = \vert det (P \, \vert _{E_{1}}) \vert
\end{displaymath}

where $P: E \to E_{2} ^{\perp}$ is the canonical projection. From the definition 
it is clear that the angle is a continuous function and that $ang (E_{1}, E_{2}) 
= 0$ iif $ E_{1} \cap E_{2} \neq \emptyset$. \\

Let $\fibvecsx$ be a 
symplectic vector bundle, where $X$ is a compact manifold and $\phi _{t}: X \to 
X$ a flow of class $C^{\infty}$, which preserves a measure $dx$. A linear symplectic cocycle is a family of 
differentiable maps $\phitast (x): \sx \to S(\phi _{t} (x))$, such that $\forall x \in X, t \in \numreal$

\begin{enumerate}
\item $\phi _{t+s} ^{\ast} = \phitast (\phi _{s} (x)) \, \phi _{s} ^{\ast} (x)$;
\item $\phi _{t} \circ \pi = \pi \circ \phi _{t}^{*}$;
\item $\forall x \in X, t \in {\mathbb R}$, $\phitast (x): \sx \to S(\phi _{t} (x))$ 
is a linear symplectic isomorphism.
\end{enumerate}

For a continuous Lagrangian distribution $\cal F$, we say that $\phitast (x)$ is $\cal F$-optical if it is in the sense of \S \ref{intro}.  \\

To prove Proposition \ref{prop-clave-prev}, we need several auxiliary lemmas. We note the Lagrange-Grassman bundle asociated to $\fibvec$ as $\Lambda (S)$. We closely follow the proof of Proposition 4.18 in Paternain \cite{pat}. \\

\begin{lema}
\label{pr-lema}
There exists $\delta > 0$, an integer $m \geq 1$ and an upper semicontinuous function
 
\begin{displaymath}
\tau : \LsSigma \times \numreal \to \{0, 1/m, 2/m, \ldots, 1 \}
\end{displaymath}

such that for all $(x, E) \in \LsSigma, t \in \numreal$, 

\begin{displaymath}
\label{prim-lema}
ang (\phi _{t + \tau} ^{\ast}(x) (E), \alpha (\phi _{t + \tau}(x) ) ) > \delta.
\end{displaymath}
\end{lema}

\demo \\
If we find $\delta > 0$ and an integer $m \geq 1$ such that if for all $(x, E, t) \in \LsSigma \times \numreal$, 

\begin{displaymath}
Q(x, E, t) = \{ i \in \entero, 0 \leq i \leq m: ang (\phi _{t + i/m} ^{\ast}(x) 
(E), \alpha (\phi _{t + i/m}(x) ) > \delta \}
\end{displaymath}

were not empty, then, clearly

\begin{displaymath}
\tau (x, E, t) = min \{i/m, i \in Q(x, E, t)\}
\end{displaymath}

has the stated properties. 

Let us assume that this is not true. Then, for every integer $m 
\geq 1$, there is a sequence $(x_{m}, E_{m}, t_{m})$ such that 

\begin{equation}
\label{eq:ang-abs}
ang (\phi _{s + t_{m}} ^{\ast} (x_{m}) (E), \alpha (\phi _{s + t_{m}}(x))) \leq 
\frac{1}{2^{m}}
\end{equation}
 
where $s \in A_{m}$, for $A_{m} = \{j / 2^{m}, j \in \entero, 0 \leq j \leq 2^{m} \}$. As a result of the compacity of $\LsSigma$, there is a subsequence converging to  $(x, E)$. Then, (\ref{eq:ang-abs}) and the continuity of the angle imply that exists $m_{k}$ such that

\begin{displaymath}
ang (\phi _{s} ^{\ast}(x)(E), \alpha (\phi_{s}(x)) ) = 0, \quad \forall s \in A_{m_{k}}
\end{displaymath}

and then, for all $s \in [0,1]$

\begin{displaymath}
\phi _{s} ^{\ast}(x)(E) \cap \alpha (\phi_{s}(x) ) \neq \{ 0 \}
\end{displaymath}

which contradicts Proposition \ref{propbp}. $\square$

\begin{lema}
\label{seg-lema}
There exists $\gamma > 0$, an integer $n \geq 1$ and an upper semicontinuous function 
\begin{displaymath}
\rho : \LsSigma \to \{0, 1/n, 2/n, \ldots, 1 \}
\end{displaymath}

such that for all $(x, E) \in \LsSigma$, 

\begin{displaymath}
ang (E, \phi _{\rho} ^{\ast} (x_{-}) \, (\alpha (x_{-}))) > \gamma
\end{displaymath}

where $x_{-} = \phi _{- \rho} (x)$.

\end{lema}

\demo \\
The proof is similar to that of Lemma \ref{prim-lema}. $\square$ \\

\begin{lema}
\label{lemlag}
For every $x \in X, t \in \numreal$, there exists a Lagrangian subspace $R_{t}(x) 
\subset S(x)$, which depends measurably on $t$ and $x$, such that

\begin{enumerate}
\item $\vert det \, \phitast (x) \vert_{R_{t}(x)} \vert = ex \, \phitast (x)$;
\item if $E \subset S(x)$, with $dim \, E = \frac{1}{2} \, dim \, 
S(x)$, then

\begin{displaymath}
\vert det \, \phitast (x) \vert_{E} \vert \geq ang (E, R_{t} ^{\perp} (x)) \, ex 
\, \phitast (x).
\end{displaymath}

\end{enumerate}
\end{lema}

\demo \\
We consider the polar decomposition

\begin{displaymath}
\phitast (x) = O_{t} (x) \, L_{t} (x)
\end{displaymath}

where $L_{t} (x): \sx \to \sx$ is symmetric and positive definite and $O_{t} (x): \sx \to S(\phi_{t}(x))$ is a linear isometry. As $\phitast (x)$ is a symplectic map, so is

\begin{displaymath}
L_{t} (x) = (\phitast (x) \, \phitast (x) ^{t}) ^{1/2}.
\end{displaymath}

As a result of this, the eigenvalues of $L_{t} (x)$ are real, and it is true that $\lambda 
_{i + n} = \lambda _{i} ^{-1}$, $i = 1, \ldots, n$. Let us assume that $\lambda_{i} \geq 1$, for $1 \leq i \leq n$. Given an almost complex structure $J(x)$ , if $v_{i}$ is an eigenvector with eigenvalue $\lambda_{i}$, then $J(x) v_{i}$ is an eigenvector with eigenvalue $\lambda _{i} 
^{-1}$, as $L_{t} \circ J = J \circ L_{t} ^{-1}$.  Then

\begin{displaymath}
\{ \vuno, \ldots, v_{n}, J(x) \vuno, \ldots, J(x) v_{n} \}
\end{displaymath}

is a base of $\sx$ and the subspace $R_{t}(x)$ generated by $\{ \vuno, \ldots, v_{n} \}$ is a Lagrangian one, on which the expansion takes place. This proves $1)$. 

To prove $2)$,  let us note that $R_{t}(x)$ and $R_{t} ^{\perp} (x)$ are invariant under $L_{t} (x)$. If $E \cap R_{t} ^{\perp} (x) = \{ 0 \}$ 
(on the contrary there is nothing to prove), then as $L_{t} 
(x) \circ P = P \circ L_{t} (x)$, with $P$ the orthogonal projection $P: S(x) \to 
R_{t}(x)$,  we have that

\begin{displaymath}
\vert det \, L_{t} (x) \vert_{R_{t}(x)} \vert  \vert det \, P \vert_{E} \vert = 
\vert det \, P \vert_{L_{t} (x)(E)} \vert  \vert det \, L_{t} (x) \vert_{E} 
\vert \leq \vert det \, L_{t} (x) \vert_{E} \vert
\end{displaymath}

which implies

\begin{displaymath}
ex \, \phitast (x) \, ang (E, R_{t} ^{\perp} (x)) \leq \vert det \, \phitast (x) 
\vert_{E} \vert 
\end{displaymath}

We now prove measurability of  $R_{t}(x)$ on $t$ and $x$. Let $\pi : \mathbb{F} \to X$ be the fibre bundle with fibres $\pi ^{-1} = \{h: S(x) \to S(x)\}$, for $h$ a linear symmetric map. For positive integers 
$p$ and $l_{i}, 1 \leq i \leq p$, we define $\mathbb{F} (l_{1}, \ldots, 
l_{p})$ as those $(x,h)$ with $p$ eigenvalues $\lambda 
_{i}$, their multiplicities being $l_{i}$. This is a Borelian set, and so is 
${\mathbb{P}} (l_{1}, \ldots, l_{p}) \subset X \times \numreal$, defined as  
${\mathbb{P}} (l_{1}, \ldots, l_{p}) = \{ (x,t), L_{t}(x) \in {\mathbb{F}} 
(l_{1}, \ldots, l_{p}) \}$. As $R_{t}(x)$ is continuous on each $\mathbb{P}$ and they are finite, the result follows. 
$\square$

\begin{lema}
\label{lema-ang}
There exists $\delta > 0$, an integer $m \geq 1$ and  measurable functions $\tau _{i} : X 
\to \{ 0, 1/m, 2/m, \ldots, 1 \}$ such that if $\tau = \tuno + \tdos$ y $x_{1} = 
\phi _{- \tuno} (x)$, $x_{2} = \phi _{\tdos} (x)$ and $\alpha ^{i} = \alpha 
(x_{i})$, then for all $x$ y $t$, 

\begin{enumerate}
\item $ang (\phi _{\tuno} ^{\ast} (x_{1}) (\alpha ^{1}), R_{t} ^{\perp} (x)) > 
\delta$;
\item $ang (\phi _{t + \tau} ^{\ast} (x_{1}) (\alpha ^{1}), \alpha ^{2}) > 
\delta$
\end{enumerate}
\end{lema}

\demo \\
If we find that, for integers $\delta _{1}, \delta_{2} > 0$, integers $m_{1}, m_{2}$ and measurable functions $\tau _{i} : X  \times \numreal \to \{ 0, 1/m_{i}, 2/m_{i}, \ldots, 1 \}$  such that $1)$ holds for $\delta = \delta _{1}$ and $2)$ for $\delta = \delta _{2}$, then the Lemma is true for $m = m_{1} \, m_{2}$ and 
$\delta = min \{\delta _{1}, \delta_{2} \}$. Taking

\begin{displaymath}
\delta _{1} = \gamma, \quad m_{1} = n, \quad \tuno (x,t) = \rho(x, R_{t} 
^{\perp} (x))
\end{displaymath} 

for $\gamma, n, \tau$ as in Lemma \ref{seg-lema} and applying it to $E = R_{t} 
^{\perp}$, we get 

\begin{displaymath}
ang (\phi _{\tuno}^{\ast} (x_{1})  (\alpha ^{1}), R_{t} ^{\perp}) > \delta_{1}
\end{displaymath}

which proves $1)$. The same idea works for the proof of $2)$, applying Lemma \ref{prim-lema} to $E = \phi _{t + \tuno + \tdos} ^{\ast} (x_{1}) (\alpha ^{1})$. $\square$

\begin{lema}
\label{lema5}
Let $\tuno, x_{1}, \alpha ^{1}$ be as in Lemma \ref{lema-ang}. Then, there exists a constant $K > 0$, such that for all $x, t$ 

\begin{displaymath}
\vert det \, \phi_{t} ^{\ast} (x_{1}) \vert_{\alpha ^{1}} \vert \geq K \, ex \,  
\phitast (x).
\end{displaymath}

\end{lema}

\demo \\
If $E = \phitast (x_{1}) (\alpha ^{1})$ we take $a > 0$ such that

\begin{displaymath}
\vert det \, \phi_{s} ^{\ast} (y) \vert _{L} \vert \geq a
\end{displaymath}

for all $s \in [0,1]$, $(y, L) \in \LsSigma$. Then, by 
$1)$ from Lemma \ref{lema-ang}

implies that

\begin{displaymath}
\vert det \, \phitast (x_{1}) \vert_{\alpha ^{1}} \vert = \vert det \, \phi _{t 
- \tuno} ^{\ast} (x) \vert_{E} \, \vert \vert det \, \phi _{\tuno} ^{\ast} 
(x_{1}) \vert_{\alpha ^{1}} \vert \geq a \, \delta \, ex \, \phi _{t - \tuno} 
^{\ast} (x). 
\end{displaymath} 

Taking $a'$ such that  $ex \, \phi _{s} ^{\ast} (x) \leq \, a'$ and for $K = a \delta a'$, we see the result is true. $\square$ \\

We now prove the main result of the section. \\

\demo \,  {\bf of Proposition \ref{prop-clave-prev} }\\
For a fixed $t$, we make the change of variables $F: X \to X$, given by $F(x) 
= x_{1} = \phi _{- \tuno(x,t)} (x)$. In each set

\begin{displaymath}
A(i) = \{ x \in X : \tuno(x,t) = i/m \}
\end{displaymath}

$F$ is injective and, as a consequence of this, for the measure $\mu = dx$ and each Borel set $S \subset A(i)$, we have $\mu (S) = \mu 
(F(S))$. If $\Phi$ is integrable in $X$, then

\begin{displaymath}
\int _{F(A(i))} \Phi \, dx = \int _{A(i)} (\Phi \circ F) \, dx.
\end{displaymath}

For $\Phi \geq 0$, then

\begin{eqnarray}
\int _{X} (\Phi \circ F) \, dx & = & \int _{\bigcup A(i)} (\Phi \circ F) \, dx = 
\sum _{i} \int _{A(i)} (\Phi \circ F) \, dx \nonumber \\ & = & \sum _{i} \int 
_{F(A(i))} \Phi \, dx \leq \sum _{i} \int _{X} \Phi \, dx = (m + 1) \int _{X} 
\Phi \, dx. \nonumber
\end{eqnarray}

Taking $\Phi = \vert det \, \phitast(x) \vert_{\alx} \vert$ and using Lemma 
\ref{lema5} 

\begin{eqnarray}
\int _{X} \vert det \, \phitast(x) \vert_{\alx} \vert \, dx & \geq & \frac{1}{m 
+ 1} \int_{X} \vert det \, \phitast (x_{1}) \vert_{\alpha ^{1}} \vert \, dx  
\nonumber \\ & \geq & \frac{K}{m + 1} \int _{X} ex \, \phitast(x) \, dx. 
\nonumber
\end{eqnarray}

which proves the Lemma for $C = \frac{K}{m + 1}$. $\square$

\section{Proof of Theorem \ref{teoa}}
\label{pruebateoa}

We state a result by Kozlovski \cite{k} and a proposition that are essential in the proofs of Theorems \ref{teoa} and \ref{teob}. 

\begin{teo} [Kozlovski \cite{k}] 
\label{Kozlovski}
For a $C^{\infty}$ diffeomorphism $f:X \to X$, where $X$ is a compact manifold

\begin{displaymath}
\htop (f) = \lim _{n \to \infty} \frac{1}{n} \log \int _{X}  ex \, (df^{n}) _{x} \, dx.
\end{displaymath}
\end{teo}

\begin{prop}
\label{hkunion}
Let $\phit: Y \to Y$ be a continuous flow on a compact manifold $Y$. Then, for $Y_{1}, Y_{2}$ closed invariants subsets in $Y$
\begin{displaymath}
\htop (\varphi, Y_{1} \cup Y_{2}) = \max _{Y_{i}, i = 1, 2} \htop (\varphi, Y_{i}).
\end{displaymath}
\end{prop}

\demo \, {\bf of Theorem \ref{teoa}} \\
To the double $D(\nepsprima)$ of $\nepsprima = H ^{-1} [e - \epsilon, e + \epsilon]$, we apply Theorem \ref{Kozlovski} and we obtain

\begin{equation}
\label{eq:expdoble}
\htop (\phit \vert _{D(\nepsprima)} ) = \lim _{t \to \infty} \frac{1}{t} \log \int _{D(\nepsprima)} \ex \, (\dphi_{t})_{x} \, dx.
\end{equation}

As

\begin{equation}
\label{eq:cortardoble}
\int _{D(\nepsprima)} \ex \, (\dphi_{t})_{x} \, dx = 2 \, \int _{\nepsprima} \ex \, (\dphi_{t})_{x} \, dx
\end{equation}

taking exponential growth rate and through (\ref{eq:expdoble}), (\ref{eq:cortardoble}) and Proposition \ref{hkunion} we get

\begin{equation}
\label{eq:expneps}
\htop (\phit \vert _{\neps} )= \lim _{t \to \infty} \frac{1}{t} \log \int _{\neps} \ex \, (\dphi_{t})_{x} \, dx
\end{equation}

where $\neps = \nepsprima \ / \,  \partial \nepsprima$. Clearly, $(d \varphi _{t})_{x}$ is an $\cal F$-optical linear symplectic cocycle, for the given distribution $\cal F$, with respecto to the flow $\phit$, for the fibre bundle $\pi: TM \vert _{\neps} \to \neps$. Proposition  \ref{prop-clave-prev} and the trivial inequality 

\begin{displaymath}
\label{eq:trivialteoa}
\ex (\dphi_{t})_{x} \geq \vert det \, (\dphi _{t})_{x} \vert _{\alx} \vert 
\end{displaymath}

turn (\ref{eq:expneps}) into

\begin{equation}
\label{eq:detdifneps}
\htop (\phit \vert _{\neps} )= \lim _{t \to \infty} \frac{1}{t} \log \int _{\neps} \vert det \, (\dphi _{t})_{x} \vert _{\alx} \vert \,dx.
\end{equation}

Let us define the function $h(\epsilon) = \htop (\phit \vert _{\neps})$. As $h$ is increasing and $h(0) = \htop (\phit \vert _{\Sigma})$, we see that

\begin{equation}
\label{eq:mayorqsigma}
\liminf _{\epsilon \to 0} h(\epsilon)  \geq \htop (\phit \vert _{\Sigma}).
\end{equation}

We use now the following result, due to Bowen \cite{bowen}.

\begin{prop} [Corolary 18, Bowen \cite{bowen}]
Let $X$ and $Y$ be compact metric spaces and  $\phit: X \to X$ a flow. If $\pi : X \to Y$ is a continuous map such that $\pi \circ \phit = \pi$, then

\begin{displaymath}
\htop (\varphi) = \sup _{y \in Y} \htop (\varphi \vert _{\pi ^{-1} (y)}). 
\end{displaymath}

\end{prop}

For $\phit \vert _{\neps}$ and $H \vert _{\neps}$, then

\begin{displaymath}
h(\epsilon) = \sup _{\epsilon _{0} \in [e - \epsilon, e + \epsilon]} \htop (\phit \vert _{H ^{-1} (\epsilon _{0})}).
\end{displaymath}

Given $r > 0$, there is an $\epsilon _{0} (\epsilon, r)$, such that

\begin{displaymath}
h(\epsilon) \leq \htop (\phit \vert _{H ^{-1} (\epsilon _{0})}) + r, \quad \epsilon _{0} \in [e - \epsilon, e + \epsilon].
\end{displaymath}

Taking upper limit when $\epsilon \to 0$ and using the fact that  topological entropy is upper semicontinuous for $C^{\infty}$ flows (Newhouse \cite{newhouse}), then

\begin{equation}
\label{eq:menorqsigma}
\limsup _{\epsilon \to 0} h(t) \leq \htop (\phit \vert _{\Sigma}).
\end{equation}

As a result of  (\ref{eq:detdifneps}), (\ref{eq:mayorqsigma}) and (\ref{eq:menorqsigma}) 

\begin{displaymath}
\htop (\phit \vert _{\Sigma}) = \lim _{\epsilon \to 0} \htop (\phit \vert _{\neps}) = \lim _{\epsilon \to 0} \lim _{t \to \infty} \frac{1}{t} \log \int _{\neps} \vert det \, (\dphi _{t})_{x} \vert _{\alx} \vert \,dx
\end{displaymath}

which proves the Theorem. $\square$

\section{Proof of Theorem \ref{teob}}
\label{pruebateob}

To prove Theorem \ref{teob}, we make descend the distribution $\cal F$ and the opticity of the Hamiltonian $H$ to a symplectic vector bundle $\fibvec$ , where the fibres are $\fibra$. \\

\begin{prop}
\label{ultimisima}
The fibre bundle given by $\fibvec$, where $S(x) = \txSigma / \xhx$, is symplectic for the form $\oms = p_{*} \omega$, where $\omega$ is the symplectic form on $M$ and $\px$ is the canonical projection. If $\altilde = \alx \cap \txSigma + \langle \xhx \rangle$, then $\altilde \subset \txSigma$ is a Lagrangian subspace of $(T_{x}M, \omega)$ and $\alsxtilde = p_{x} (\altilde)$ is a Lagrangian subspace of $(S(x), \omega _{S})$.
\end{prop}

\demo \\

For $v \in \txSigma$, we note $p_{x} (v) = [v]$. Then, for $\vuno, \vdos \in \txSigma$, we define

\begin{displaymath}
\oms ([\vuno], [\vdos]) = \omSigma (\vuno, \vdos).
\end{displaymath}
As $\omega _{\Sigma}$ degenerates on $\xhx$ it is clear that $\omega _{S}$ is a symplectic form on $S(x)$. The subspace $\altilde \subset \txSigma$ annihilates $\omega$, as for $v' \ in \altilde$, then $v' = v \, + \, a \langle \xhx \rangle$, with $v \in \alx$ and $a \in \numreal$. In case $\alx \subset \txSigma$ and $\langle \xhx \rangle \not \in \alx$, then  $dim  \, \altilde \, = \, n \, + \, 1$. But this implies that $\altilde$ is an $(n+1)$-dimensional Lagrangian subspace, as $\altilde \, = \, (\altilde) ^{\perp}$. So $\altilde$ has dimension $n$ and this proves it is a Lagrangian subspace. $\square$

\begin{prop}
\label{bajaopt}
Let $\fstilde$ be the Lagrangian distribution induced in $\fibvec$ by  Proposition \ref{ultimisima}. If $H \in C^{\infty}(M, \numreal)$ is an $\cal F$-optical Hamiltonian, then $\widetilde{\dphi _{t}} (x)$ is an $\fstilde$-optical linear symplectic cocycle.
\end{prop}

\demo \\
By the definition of the optical property, it follows that for a curve

\begin{displaymath}
\lambda (t) = S(t) \alpha (\varphi _{t} (x)), \qquad S(0) = Id, \, S(t) \in Sp \, (T _{\varphi_{t}(x)}M, \omega)
\end{displaymath}

in the fibre bundle $\Theta: \cal A \to M$ and to the vector $\dot{\lambda}(0) \in T_{\alx} \cal A$ given by 

\begin{displaymath}
\dot{\lambda}(0) = \frac{d}{dt} \, (\phit)_{*} \, (\alpha (\varphi _{-t} (x)) \, \vert _{t = 0}
\end{displaymath}

we can associate a bilinear symmetric positive definite form

\begin{equation}
\label{eq:bil2}
(\zetaeta) \mapsto \omega (\zeta, \dot{S}(0) \eta).
\end{equation}

Each $\altilde$ contains $\xhx$ and is included in $\txSigma$, which implies that the form (\ref{eq:bil2}) is positive semidefinite, as $\omega \vert_{\Sigma}$ degenerates on the field. Descending to the quotient $\sx$,  $\oms$ from  Proposition \ref{ultimisima} becomes symplectic and then

\begin{displaymath}
([\zeta], [\eta]) \mapsto \omega([\zeta], [\dot{S}(0) \eta]), \qquad [\zeta], [\eta] \in \alsxtilde
\end{displaymath}

is the symmetric bilinear form which proves the Proposition. $\square$ \\

We note by $ex _{\txSigma}$ the expansion on $\txSigma$ and by $\ex _{S(x)}$ the expansion on $S(x)$. \\

The following lemmas relate expansions and determinants of the differential of the flow on  $\txSigma$ and $S(x)$. We only prove the second one, as the proof of the first is easy.

\begin{lema}
\label{primerlemafinal}
Let $\pi : E \to X$ be a fibre bundle on a compact manifold $X$. If $g_{1}$ and $g_{2}$ are two continuous Riemannian metrics on $E$, then there is a constant  $K > 0$, such that for the linear map $L: E(x) \to E(y)$ 

\begin{displaymath}
\ex _{g_{1}} L \leq K \, \ex _{g_{2}} L.
\end{displaymath}

\end{lema}

\begin{lema}
\label{segundolemafinal}
For a symplectic manifold $(M, \omega)$, endowed with a Lagrangian distribution   $\cal F$ and an  $\cal F$-optical Hamiltonian flow $\phit$ that admits a continuous invariant distribution of hyperplanes $\xhx$ on $\txSigma$, then there exists constants $K_{1}, K_{2} > 0$ such that 

\begin{eqnarray*}
\ex _{\txSigma} \, \dphi _{t} (x) \leq K_{1} \, \ex _{S(x)} \, \widetilde{\dphi _{t}} (x) \\
\vert det \, \widetilde{\dphi _{t}} (x) \vert _{\alsxtilde} \vert \leq K_{2} \, \vert det \, \dphi _{t} (x) \vert _{\altilde} \vert.
\end{eqnarray*}

\end{lema}

\demo  \\
The symplectic form $\omega$ is nondegenerate on $T(x)$ due to the transversality with respect to $\xhx$. Then, there exists and almost complex structure $J$ that induces an inner product $\langle .,. \rangle _{T}$ in $T$. We extend it to $\txSigma$ through

\begin{displaymath}
\label{prodint}
\langle v,w \rangle _{\txSigma} = \left \{ \begin{array} {c@{\quad \, si \,\quad}l}
\langle v,w \rangle _{T} & v,w \in T \\ 0 & v \in \langle \xhx \rangle, w \in T(x) \\
1 & v = w = \xhx
\end{array} \right.
\end{displaymath}

Invariance under the flow allows us to choose a base $\{ \xhx, t_{1}, \ldots, t_{2n-2} \}$ for $\txSigma$, where $\{ t_{i} \} _{1 \leq i \leq 2n-2}$ is a base of $T(x)$, such that 

\begin{displaymath}
\label{eq:matriz}
(d\phit)_{x} \vert _{\Sigma} = \left( \begin{array}{cc} 1 & 0 \\ 0 & (d\phit)_{x} \vert _{T(x)} \end{array} \right).
\end{displaymath}

Then

\begin{displaymath}
\label{eq:igexpSigmaT}
\ex _{\txSigma} (d\phit) _{x} = \ex _{T(x)} (d\phit) _{x}.
\end{displaymath}

as the expansion is the maximum of the determinant on all minors of any dimension on this matrix.

We introduce an inner product on $S(x)$, such that the projection $P \vert _{T(x)}$ is an isometry. It is clear that if $g_{1}, g_{2}$ are the metrics given by the inner products on $T(x)$ and $S(x)$, then

\begin{displaymath}
\ex _{T(x)} ^{g_{1}} (\dphi _{t}) _{x} = \ex _{S(x)} ^{g_{2}} (\widetilde{\dphi _{t}}) _{x}.
\end{displaymath}

Lemma \ref{primerlemafinal} leads us to the first inequality we wanted to prove. The proof of the second one, leans on a similar method. $\square$ \\

\demo \, {\bf of Theorem \ref{teob}} \\
Proposition \ref{bajaopt} implies that $(\widetilde{d\phit})_{x}$ is an $\fstilde$-optical linear symplectic cocycle, for the distribution given by  Proposition \ref{ultimisima}, with respect to the flow $\phit \vert _{\Sigma}$, on the symplectic vector bundle $\fibvec$. Proposition \ref{prop-clave-prev} and the trivial inequality 

\begin{displaymath}
\ex (\widetilde{d\phit})_{x} \geq \vert det \, (\widetilde{d\phit})_{x} \vert _{\alsxtilde} \vert
\end{displaymath}

lead us to 

\begin{equation}
\label{eq:igualintexpdif}
\liminf _{t \to \infty} \frac{1}{t} \log \int _{\Sigma} \ex _{S(x)} \, (\widetilde{d\phit})_{x} \, dx = \liminf _{t \to \infty} \frac{1}{t} \log \int _{\Sigma} \vert det \, (\widetilde{d\phit})_{x} \vert _{\alsxtilde} \vert \,dx. 
\end{equation}

Then, as a result of Kozlovski's Theorem \ref{Kozlovski}, the first inequality from Lemma \ref{segundolemafinal}, (\ref{eq:igualintexpdif}) and the second inequality form Lemma \ref{segundolemafinal}, 

\begin{eqnarray}
\label{eq:primerdes}
\htop (\phit \vert _{\Sigma}) & = & \liminf _{t \to \infty} \frac{1}{t} \log \int _{\Sigma} \ex _{\txSigma} (\dphi_{t})_{x} \, dx  \nonumber \\ & \leq &  \liminf _{t \to \infty} \frac{1}{t} \log \int _{\Sigma} \ex _{S(x)} \, (\widetilde{\dphi _{t}}) _{x} \, dx \nonumber \\ & = & \liminf _{t \to \infty} \frac{1}{t} \log \int _{\Sigma} \vert det \, (\widetilde{\dphi _{t}}) _{x} \vert _{\alsxtilde} \vert \,dx \nonumber \\ & \leq & \liminf _{t \to \infty} \frac{1}{t} \log \int _{\Sigma} \vert det \, (\dphi _{t})_{x} \vert _{\altilde} \vert \,dx.
\end{eqnarray}

Similar methods imply that

\begin{eqnarray}
\label{eq:segundades}
\htop (\phit \vert _{\Sigma}) & = & \lim _{t \to \infty} \frac{1}{t} \log \int _{\Sigma} \ex _{\txSigma} (\dphi_{t})_{x} \, dx  \nonumber \\ & \geq & \limsup _{t \to \infty} \frac{1}{t} \log \int _{\Sigma} \vert det \, (\dphi _{t})_{x} \vert _{\altilde} \vert \,dx
\end{eqnarray}

which through (\ref{eq:primerdes}) and (\ref{eq:segundades}) say that 

\begin{displaymath}
\htop (\phit \vert _{\Sigma}) = \lim _{t \to \infty} \frac{1}{t} \log \int _{\Sigma} \vert det \, (\dphi _{t})_{x} \vert _{\altilde} \vert \,dx
\end{displaymath}

which proves the Theorem. $\square$

\end{document}